\font\sets=msbm10.
\font\script=eusm10.
\def\1{{\bf 1}}
\def\square{\hbox{\vrule\vbox{\hrule\phantom{s}\hrule}\vrule}}
\def\sgn{{\rm sgn}}
\def\defineq{\buildrel{def}\over{=}}
\def\C{\hbox{\sets C}}
\def\N{\hbox{\sets N}}

\def\R{\hbox{\sets R}}
\def\Z{\hbox{\sets Z}}
\def\divisor{\hbox{\bf d}}
\def\sporadic{\hbox{\script S}}
\def\exceptional{\hbox{\script E}}
\def\DFTf{\widehat{f_{_{\cal N}}}}

\def\gQ{{g_{_Q}}}

\par
\centerline{\bf Sieve functions in arithmetic bands, II}
\bigskip
\centerline{Giovanni Coppola and Maurizio Laporta}\footnote{}{{\bf Math Subject Classification 2010 :} 11N37 \quad {\bf Keywords:} arithmetic progressions, short intervals}

\bigskip	
\bigskip
\bigskip

\par
\noindent {\bf Abstract.} An arithmetic function $f$ is called a {\it sieve function} of {\it range} $Q$ if its Eratosthenes transform $g=f\ast\mu$ is supported in $[1,Q]\cap\N$, where $g(q)\ll_{\varepsilon} q^{\varepsilon}$ ($\forall\varepsilon>0$). We continue our study of the distribution of such functions over short {\it arithmetic bands}, $n\equiv ar+b\, (\bmod\,q)$, with $1\le a\le H=o(N)$ and $r,b\in\Z$ such that ${\rm g.c.d.}(r,q)=1$. In particular, we discuss the optimality of some results. 

\bigskip
\bigskip

\par
\centerline{\bf 1. Introduction.}
\smallskip
\par
\noindent
Given an arithmetic function $f:\N \rightarrow \C$, for every
$N\in\N$ and  $\alpha \in \R$ let us set
$$
\DFTf(\alpha)\defineq \sum_{n\sim N}f(n)e(n\alpha), 
$$
\par
\noindent
where  $n\sim N$ means that $n\in{\cal N}\defineq(N,2N]\cap \N$ and $e(\alpha)\defineq e^{2\pi i\alpha}$, as usual.
If $f$ is the convolution product of $g$ and 
the constantly $1$ function, i.e. 
$$
f(n)=(g\ast \1)(n)=\sum_{d|n}g(d),
$$
\par
\noindent
we say, with Wintner [W], that $g=f\ast\mu$ is the {\it Eratosthenes transform} of $f$ (where $\mu$ is the well-known M\"{o}bius function). 
We call $f$ a {\it sieve function} of {\it range} $Q$ if its Eratosthenes transform
$g$ is {\it essentially bounded}, namely $g(d)\ll_{\varepsilon} d^{\varepsilon}$, $\forall\varepsilon>0$, and vanishes outside $[1,Q]$ for some $Q\in\N$, that is to say,
$$
f(n)=\sum_{{d|n}\atop {d\le Q}}g(d). 
$$
\par
\noindent
Note that $f=g\ast \1$ is essentially bounded 
if and only if so is $g$.

As usual, $\ll$ is Vinogradov's notation, synonymous to Landau's $O$-notation. In particular, $\ll_{\varepsilon}$ means that the implicit constant might 
depend on an arbitrarily small and positive real number $\varepsilon$, which might change at each occurrence. We also write $\gQ\defineq g\cdot \1_Q$ to mean that $g$ vanishes outside $[1,Q]$, the function $\1_Q$ being the indicator of $[1,Q]\cap\N$.
In the above notation, we have set $f_{_{\cal N}}\defineq f_{2N}-f_N$.

In [CL2] we have established a formula that relates
the so-called $\ell$th {\it Ramanujan coefficient} of a real sieve function $f$ of range $Q$, i.e.
$$
R_{\ell}(f)\defineq \sum_{d\equiv 0\, (\ell)}{{\gQ(d)}\over d},
$$
\par
\noindent
to the values of $\DFTf$ attained at rational numbers; 
hereafter we write $n\equiv a\, (q)$ to abbreviate $n\equiv a\, (\bmod\, q)$. More precisely, Lemma 3.1 [CL2] states that
$$
\DFTf(j/\ell)=R_{\ell}(f)N+O_{\varepsilon}((\ell Q)^{\varepsilon}(Q+\ell)), 
\quad \enspace \forall \ell>1, \forall j\in\Z_{\ell}^{*}, 
\leqno{(1)}
$$
\par
\noindent
which holds uniformly in the complete set  $\Z_{\ell}^{*}$ of reduced residue classes modulo $\ell$. Moreover, it easily seen that (compare also [C])
$$
R_{\ell}(f)
={1\over {\ell}}\sum_{m\le Q/\ell}{{g(\ell m)}\over m}
\ll_{\varepsilon} {{Q^{\varepsilon}}\over {\ell}}.
\leqno{(2)}
$$

\par				
Here we continue our study of the distribution of a real sieve funtion $f$ over short {\it arithmetic bands}, i.e. 
$$
\bigcup_{1\le a\le H}\{n\in(N,2N]: n\equiv a\, (\bmod\,q)\},\ \hbox{with}\ H=o(N).
$$
\par
\noindent
Let us recall that in [CL2] we have proved that the inequality (hereafter, we omit $a\ge 1$ in sums like $\sum_{a\le H}$)
$$
T_f(q,N,H)\defineq\sum_{a\le H}\sum_{{n\sim N}\atop {n\equiv a\, (q)}}f(n)
-{H\over q}\DFTf(0)\ll_{\varepsilon}N^{\varepsilon}(N/q+q+Q)
\leqno{(3)}
$$
\par
\noindent
holds for every real sieve function $f$ of range $Q\ll N$, after assuming that 
$H=o(q)$, as $q\to \infty$, and $q=o(N)$, as $N\to \infty$. Such conditions are required in order to
avoid overlapping and {\it sporadicity} of the arithmetic bands, respectively. 
By a straightforward application of $(1)$ and $(2)$ here we generalize the previous inequality for
$$
T_f(q,r,b,N,H)\defineq \sum_{a\le H}\sum_{{n\sim N}\atop {n\equiv ar+b\, (q)}}f(n)-{H\over q}\DFTf(0), 
$$
\par
\noindent
where $r,b\in\Z$ are such that $(r,q)=1$ (hereafter $(r,q)={\rm g.c.d.}(r,q)$, as usual). 
In particular, note that $T_f(q,1,0,N,H)=T_f(q,N,H)$.
\smallskip
\par
\noindent {\bf Theorem.} {\it Let } $q,N,H,Q\in\N$ {\it be such that } $Q\ll N$, $H=o(q)$ {\it and } $q=o(N)$. {\it For every sieve function } $f:\N \rightarrow \R$ {\it of range } $Q$ {\it and all } $r,b\in\Z$ {\it such that } $(r,q)=1$, {\it one has} 
$$
T_f(q,r,b,N,H)\ll_{\varepsilon}N^{\varepsilon}(N/q+q+Q). 
\leqno{(4)}
$$
\par
\noindent {\bf Proof.} By the orthogonality of additive characters
$e_q(t)\defineq e(t/q)$,\ ($q\in \N$,  $t\in \Z$), we get 
$$
T_f(q,r,b,N,H)={1\over q}\sum_{a\le H}\sum_{n\sim N}f(n)\sum_{k\le q}e_q(k(ar+b-n))-{H\over q}\DFTf(0) 
= {1\over q}\sum_{k<q}\sum_{a\le H}e_q(k(ar+b))\DFTf(-k/q)
$$
$$
={1\over q}\sum_{{\ell|q}\atop {\ell>1}}\sum_{j\in \Z_{\ell}^{*}}\DFTf(-j/\ell)\sum_{a\le H}e_{\ell}(j(ar+b)),
$$
\par
\noindent
where we have set $\ell=q/(k,q)$ and $j=k/(k,q)$. Since $(r,q)=1$, for any $\ell|q$ there exists an integer $\overline{r}$ such that $r\overline{r}\equiv 1\, (\bmod\, \ell)$. Therefore we write
$$
T_f(q,r,b,N,H)={1\over q}\sum_{{\ell|q}\atop {\ell>1}}\sum_{j\in \Z_{\ell}^{*}}
\DFTf\left(-{j\overline{r}\over {\ell}}\right)
e_{\ell}(j\overline{r}b)\sum_{a\le H}e_{\ell}(ja).
$$
\par
\noindent
By applying $(1)$, $(2)$ and the well-known inequality (see [Da], Ch.25)
$$
\sum_{V_1<v\le V_2}e(v\alpha)\ll \min\Big(V_2-V_1,{1\over {\Vert\alpha\Vert}}\Big),
$$
\par
\noindent
we conclude that
$$
T_f(q,r,b,N,H)\ll_{\varepsilon} {1\over q}\sum_{\ell|q,\ell>1}\Big(|R_{\ell}(f)|N+(\ell Q)^{\varepsilon}(Q+\ell)\Big)\sum_{j\in \Z_{\ell}^{*}}{1\over {\Vert j/\ell\Vert}}
$$
$$
\ll_{\varepsilon} {{N^{\varepsilon}}\over q}\sum_{\ell|q,\ell>1}\Big({N\over {\ell}}+Q+\ell\Big)\ell 
\ll_{\varepsilon} N^{\varepsilon}\Big({N\over q}+Q+q\Big),
$$
\par
\noindent
that is $(4)$. \hfill $\square$ 

\medskip

\par				
\noindent {\bf Remark 1.} Besides $(4)$, from the previous proof it transpires that also the upper bound of
$$
T_f(q,r,b,N,H)={1\over q}\sum_{{\ell|q}\atop {\ell>1}}\sum_{j\in \Z_{\ell}^{*}}\DFTf\left(-{j\overline{r}\over {\ell}}\right)e_{\ell}(j\overline{r}b)\sum_{a\le H}e_{\ell}(ja) 
\ll {1\over q}\sum_{{\ell|q}\atop {\ell>1}}\ell \sum_{j\le {\ell/2}\atop {(j,\ell)=1}}{1\over j}\max_{j\in \Z_{\ell}^*}|\DFTf(j/\ell)|
$$
\par
\noindent
is independent of both $b$ and $r$ such that $(r,q)=1$ (that is to say, it is the same upper bound obtained for $r=1$ and $b=0$). 

\medskip

\par
\noindent {\bf Remark 2.} Recalling that here trivial bound means $\ll N^{1+\varepsilon}H/q$, both $(3)$ and $(4)$ are non-trivial once the {\it width} $\theta\defineq {{\log H}\over {\log N}}$ is positive, $q\ll \sqrt{N^{1-\delta}H}$ and $qQ\ll N^{1-\delta}H$, for a suitable \enspace $\delta>0$. Since it is assumed that $Q\ll q\ll Q$ hereafter, we get a bound $Q\ll \sqrt{N^{1-\delta}H}$ and thus go beyond the square-root of $N$ (for $\theta>0$). Consistently with the terminology introduced in [CL2], we stop at the {\it barrier} ${{\log Q}\over {\log N}}<{{1+\theta}\over 2}$. 

\bigskip

In the final section of [CL2] we compared our study of the distribution of sieve functions in arithmetic bands to the classical results on the distribution in  arithmetic progressions, which deal mostly with the overcoming of the so-called {\it level} $1/2$. 
In fact, our present results, insofar they generalize our previous ones, already go beyond level $1/2$. 
\par
Here we discuss the possibility of going beyond $1/2+\theta/2$ in the present contest of the sieve functions in arithmetic bands, namely by taking $Q>\sqrt{N^{1+\delta}H}$ for a certain small and fixed $\delta>0$, so that $N^{1+\varepsilon}H/Q=o(Q)$, provided 
$\delta>\varepsilon$. 
Indeed, in section 3 we exhibit a particular sieve function whose behavior in arithmetic bands confirms the optimality of such level. 

\bigskip
\bigskip

\par
\centerline{\bf 2. Length inertia in arithmetic bands.}
\smallskip
\par
\noindent
In [CL1] we studied the so-called {\it length inertia} property for 
{\it weighted Selberg integrals} (see [CL2] for the link between such integrals and the distribution of a sieve function in arithmetic bands). Such a property permits transfer of
non-trivial bounds in short intervals of length $h$, say, to 
similar bounds in longer intervals of length $H=\infty(h)$ (that is $h=o(H)$, as $H\to \infty$). Here we show that a length inertia property holds also for the distribution of a sieve function in arithmetic bands. Indeed, we have
$$
T_f(q,N,[H/h]h)\defineq \sum_{a\le [H/h]h}\sum_{{n\sim N}\atop {n\equiv a\, (q)}}f(n)-{[H/h]h\over q}\DFTf(0)
=\sum_{j\le [H/h]}\Big(\sum_{(j-1)h<a\le jh}\sum_{{n\sim N}\atop {n\equiv a\, (q)}}f(n)-{h\over q}\DFTf(0)\Big)
$$
$$
=\sum_{j\le [H/h]}\Big(\sum_{a\le h}\sum_{{n\sim N}\atop {n\equiv a+(j-1)h\, (q)}}f(n)-{h\over q}\DFTf(0)\Big)
=\sum_{j\le [H/h]}T_f(q,1,(j-1)h,N,h).
$$

\bigskip
\bigskip

\par
\centerline{\bf 3. Optimality of the level.}
\smallskip
\par
\noindent
Let us set \thinspace $\sgn(x)\defineq x/|x|$ \thinspace for all \thinspace $x\in\R\setminus\{0\}$ \thinspace and \thinspace $\sgn(0)\defineq 0$.
Then, for any fixed $q\in \N\cap(Q,2Q]$, we define the arithmetic function $g$  as
$$
g(d)=g(d,q,N,H)\defineq \sgn\Big(\sum_{a\le H}\Big(\sum_{{m\sim N/d}\atop {md\equiv a\, (q)}}1-{1\over {q}}\sum_{m\sim N/d}1\Big)\Big), 
$$
\par
\noindent
if \thinspace $d\in\N\cap(Q,2Q]$ \thinspace and \thinspace $g(d)\defineq 0$ \thinspace otherwise. It is plain that $g$ is the Eratosthenes transform of the sieve function \thinspace $f(n)=f(n,q,N,H)\defineq \sum_{d|n}g(d)$ \thinspace of range $2Q$. By noting that \thinspace $\sgn(x)x=|x|$ \thinspace for all \thinspace $x\in\R$, we write
$$
|T_f(q,N,H)|=\Big|\sum_{a\le H}\Big(\sum_{{n\sim N}\atop {n\equiv a\, (q)}}f(n)-{1\over {q}}\sum_{n\sim N}f(n)\Big)\Big|
=\Big|\sum_{d\sim Q}g(d)\sum_{a\le H}\Big(\sum_{{m\sim N/d}\atop {md\equiv a\, (q)}}1-{1\over {q}}\sum_{m\sim N/d}1\Big)\Big|
$$
$$				
=\sum_{d\sim Q}\Big|\sum_{a\le H}\Big(\sum_{{m\sim N/d}\atop {md\equiv a\, (q)}}1-{1\over {q}}\sum_{m\sim N/d}1\Big)\Big|.
$$
\par
\noindent
In order to show that \thinspace $|T_f(q,N,H)|\gg NH/q$, we set 
$$
\sporadic=\sporadic(q,Q,H,N)\defineq\Big\{ d\in \N\cap(Q,2Q] : \sum_{a\le H}\sum_{{m\sim N/d}\atop {md\equiv a\, (q)}}1\ge 1\Big\},
$$
$$
\exceptional\defineq \left(\N\cap(Q,2Q]\right)\setminus \sporadic
$$
\par
\noindent
and prove that \thinspace $|\sporadic|=o(Q)$, which in turns yields \thinspace $|\exceptional|\gg Q$. Indeed, from the latter inequality we observe that 
$$
|T_f(q,N,H)|=\sum_{d\in\sporadic} \Big|\sum_{a\le H} \Big(\sum_{{m\sim N/d}\atop {md\equiv a\, (q)}}1-{1\over {q}} \sum_{m\sim N/d}1\Big)\Big| 
             + \sum_{d\in \exceptional}{H\over {q}}\sum_{m\sim N/d}1
$$
$$
\gg |\exceptional|{{NH}\over {qQ}}
\gg {{NH}\over {q}}.
$$
\par
\noindent
Now let us prove that \thinspace $|\sporadic|=o(Q)$. To this end, after recalling that the divisor function \thinspace $\divisor(n)\defineq \sum_{d|n}1$ \thinspace is essentially bounded and \thinspace $q>Q$, we note that 
$$
|\sporadic|\le \sum_{d\in \sporadic_{q}}\sum_{a\le H}\sum_{{m\sim N/d}\atop {md\equiv a\, (q)}}1
\le \sum_{a\le H}\sum_{{n\sim N}\atop {n\equiv a\, (q)}}\divisor(n)
\ll_{\varepsilon} {{N^{1+\varepsilon}H}\over Q}.
$$
\par
\noindent
If \thinspace $Q>\sqrt{N^{1+\delta}H}$ \thinspace for a certain \thinspace $\delta>0$, then \thinspace $N^{1+\varepsilon}H/Q=o(Q)$ \thinspace once \thinspace $\varepsilon<\delta$, that yields the desired conclusion. 

\bigskip
\bigskip

\par
\centerline{\bf 4. Concluding remarks.}
\smallskip
\par
\noindent
Sieve functions are ubiquitous in analytic number theory. Besides the truncated divisor sum $\Lambda_R$ (see [G]), that is a linear combination of sieve functions of range $R$ (see [CL2] for our remarks on $\Lambda_R$), we quote the so-called {\it restricted} divisor function 
$$
\tau_Q(n)\defineq(\1_Q\ast\1)(n)=\sum_{d|n,d\le Q}1,
$$
\par
\noindent
whose Eratosthenes transform is the indicator $\1_Q$ of $[1,Q]\cap\N$.
We refer the reader to [T] for results on the distribution of $\tau_Q$ in short arithmetic progressions. Here we wish to stress that in [T] one finds 
also conjectures and average results concerning the distribution in arithmetic bands
of the more complicated function  
$$
\tau_{Q,R}(n)\defineq(\1_Q\ast\1_R)(n)=\sum_{{(d,t)}\atop {dt=n}}\1_Q(d)\1_R(t).
$$
\par
\noindent
Such an essentially bounded function is linked to the pair correlation problem for fractional parts of the quadratic function
$\alpha k^2$, with $k\in\N$ and $\alpha\in\R$
(compare also [S]). While in [T], Conjecture 1.2, it is pursued the research of an upper bound for
$$
\sum_{a\le H}\sum_{n\equiv ar\, (q)}\tau_{Q,R}(n)
-{HQR\over q},\ \hbox{with}\ (r,q)=1,
$$
\par
\noindent
under suitable conditions on $H,Q,R$, here our Theorem leads to an asymptotic formula for the inverse Eratosthenes transform of $\tau_{Q,R}$ in arithmetic bands, namely for
$$
\sum_{a\le H}\sum_{n\equiv ar\, (q)}(\tau_{Q,R}\ast\mu)(n). 
$$

\par				
\noindent
In [CL3] we established an asymptotic inequality for the exponential sum associated to the {\it localized} divisor functions, a family  of functions including the aforementioned $\tau_Q$, $\tau_{Q,R}$, and the standard divisor function
$d_k$, $k\ge 2$ (recall that
$d_k(n)$ is the number of ways to write $n$ as a product of $k$ positive integers).
The particular instance of such an inequality for $d_k$ is
$$
\sum_{n\sim N}d_k(n)e(n\alpha)\ll_{k,\varepsilon} (Nq)^{\varepsilon}(N/q+q+N^{1-1/k}),
$$
\par
\noindent
for all $\alpha\in [a/q-1/q^2,a/q+1/q^2]$, with $(a,q)=1, q>1$.
Somehow, this can be regarded as being analogous to the inequality which follows 
by combining $(1)$ with $(2)$. Such a circumstance is remarkable even in view of the fact that any
$d_k$ falls short of being a sieve function, with a sort of Eratosthenes transform which turns out to be the sum of $d_{k-1}$ plus some restricted divisor functions (see the last section of [CL1]).

Finally, since the function proposed in section 3 seems to be rather artificial in that it depends on a fixed modulus $q\in \N\cap(Q,2Q]$, an intriguing open question to ask is how many standard sieve functions might support the optimality of the level accomplished by our Theorem and the results of [CL2].

\bigskip
\bigskip

\par
\centerline{\bf References}
\medskip
\item{[C]} G. Coppola, {\sl On some lower bounds of some symmetry integrals}, Afr. Mat. {\bf 25}, issue 1 (2014), 183--195. $\underline{\tt MR\thinspace \!:\!3165958}$ 
\smallskip
\item{[CL1]} G. Coppola and M. Laporta, {\sl Symmetry and short interval mean-squares}, (2016), arXiv: 1312.5701 (submitted)
\smallskip
\item{[CL2]} G. Coppola and M. Laporta, {\sl Sieve functions in arithmetic bands}, Hardy-Ramanujan J. {\bf 39}, (2016), 21-37
\smallskip
\item{[CL3]} G. Coppola and M. Laporta, {\sl A note on the exponential sums of the localized divisor functions}, to appear in PALANGA 2016 Conference Proceedings 
\smallskip
\item{[Da]} H. Davenport, {\sl Multiplicative Number Theory}. 3rd edition, GTM {\bf 74}, Springer, New York, 2000.
\smallskip
\item{[G]} D.A. Goldston, {\sl On Bombieri and Davenport's theorem concerning small gaps between primes}, Mathematika, {\bf 39} (1992), 10--17
\smallskip
\item{[S]} I.E. Shparlinski, {\sl On the restricted divisor function in arithmetic progressions}, Rev. Mat. Iberoam., {\bf 28}, (2012), 231-238
\smallskip
\item{[T]} J.L. Truelsen, {\sl Divisor problems and the pair correlation for the fractional parts of $n^2\alpha$}, Int. Math. Res. Not. IMRN, 2010 (2010), 3144-3183
\smallskip
\item{[W]} A. Wintner, {\sl Eratosthenian Averages}, Waverly Press, Baltimore, MD, 1943 

\bigskip
\bigskip
\bigskip
	
\par
\noindent
Giovanni Coppola
\par
\noindent
Home address \negthinspace : \negthinspace Via Partenio \negthinspace 12
\par
\noindent
83100, Avellino, ITALY 
\par
\noindent
{\it e-page:} {\tt www.giovannicoppola.name} 
\par
\noindent
{\it e-mail:} {\tt giocop@interfree.it}

\medskip

\par
\noindent
Maurizio Laporta
\par
\noindent
Universit\`a degli Studi di Napoli "Federico II",
\par
\noindent
Dipartimento di Matematica e Applicazioni "R.Caccioppoli",
\par
\noindent
Complesso di Monte S.Angelo, Via Cinthia - 80126, Napoli, ITALY
\par
\noindent
{\it e-mail:} {\tt mlaporta@unina.it}

\bye